\documentclass[11pt]{article}

\usepackage{amsmath, amssymb} 
\usepackage{mathptmx} 
\usepackage{authblk}
\usepackage[a4paper, margin=1in]{geometry}

\def\dj{d\kern-0.3em\char"16\kern-0.08em}
\def\Dj{\mbox{\raise0.3ex\hbox{-}\kern-0.3em D}}

\usepackage[bottom]{footmisc}
\usepackage{theorem}

\newtheorem{theorem}{Theorem}[section]
\newtheorem{lemma}[theorem]{Lemma}

\theoremstyle{definition}

\theoremstyle{remark}

\newcommand{\dbar}{{\mkern3mu\mathchar'26\mkern-12mu d}}

\def\R{\mathbb R} \def\N{\mathbb N}   
\def\vp{\varphi}

\def\be{\begin{equation}}
\def\ee{\end{equation}}

\providecommand{\keywords}[1]{\textbf{\textit{Keywords:}} #1}   

\providecommand{\subjclass}[1]{\textbf{\textit{Mathematics Subject Classification:}} #1}

\begin{document}

\title{Continuity of pseudodifferential operators with nonsmooth\\ symbols on mixed-norm Lebesgue spaces}

\author[1]{Ivan Ivec}

\affil[1]{Faculty of Metallurgy, University of Zagreb, Aleja narodnih heroja 3, Sisak, Croatia\\

iivec@simet.unizg.hr}

\maketitle

\begin{abstract}
Mixed-norm Lebesgue spaces found their place in the study of some questions in the theory of partial differential equations,
as can be seen from recent interest in the continuity of certain classes of pseudodifferential operators
on these spaces. In this paper, we use some recent advances in the pseudodifferential calculus for nonsmooth symbols to
prove the boundedness of pseudodifferential operators with such symbols on mixed-norm Lebesgue spaces.
\end{abstract}



\let\thefootnote\relax\footnote{This work has been supported in part by the Croatian Science Foundation
under the project 2449 MiTPDE.}

\section{Introduction}

Continuity of pseudodifferential operators with symbols in the nonsmooth H\"ormader class
${\rm S}^m_{\rho,\delta,N,N'}$ (the definition is recalled at the beginning of Section 3) on mixed-norm Lebesgue spaces is our main interest in this paper.
Regarding the continuity on classical Lebesgue spaces, there is a famous result by Coifman and Meyer \cite[Theorem 7]{cm}, namely
that for $0\leq\delta\leq\rho\leq1$, $\delta<1$ and $m=0$ it is enough to have $N,N'>\frac d2$ to obtain the continuity on $L^2(\R^d)$.
Moreover, in H\"ormander \cite{HorP} it was shown that we have the following necessary condition
for continuity on $L^p(\R^d)$ spaces:
\begin{equation}\label{lpcond}
m\leq-d(1-\rho)\Big|\frac12-\frac1p\Big|\,.
\end{equation}
From this condition, it is clear that for the zeroth-order operators ($m=0$) and $p\not =2$ we can obtain the continuity result only in the case $\rho=1$. 
An interesting result for nonsmooth symbols in the $L^p$ setting is obtained in \cite[Theorem 2]{apv}, where the norm of an operator is estimated in terms of the norm of its symbol.
That was useful for construction of certain variants of H-distributions (introduced in \cite{am}), both in \cite{apv} and \cite{iv}. 
The condition (\ref{lpcond}) proved to be relevant also for continuity on Lebesgue spaces on compact and graded Lie groups \cite{R1, R2}.
We investigate this further for an arbitrary mixed-norm Lebesgue space $L^{\bf p}(\R^d)$, ${\bf p}\in \langle1,\infty\rangle^d$. Mixed-norm spaces are suitable to describe properties of functions that have different growth with respect to different variables. They are naturally connected with problems for estimating solutions of partial differential equations modeling
physical processes depending both on space and time. For instance, they were used in the famous Strichartz estimates for solutions of the Schr\"{o}dinger equation \cite{CZ}. The definition and some properties of mixed-norm Lebesgue spaces are given in \cite{BePa} and recalled in the Appendix of \cite{AIV}.

In \cite{AI, AIV} the boundedness of pseudodifferential operators on mixed-norm Lebes\-gue spaces 
have already been obtained, but these results cover only the case of infinitely smooth symbols together with $m=0$. Also, 
a detailed discussion on similar results from the literature was given there.

In the next section we recall some already known results that we need as ingredients in our proof:
sufficient conditions for boundedness of a linear operator $A:L^\infty_c(\R^d)\to L^1_{loc}(\R^d)$ on mixed-norm Lebesgue spaces, obtained 
in \cite{AIV}, and pseudodifferential calculus for the adjoint of the pseudodifferential operator with a nonsmooth symbol, obtained in \cite{IV}.
In the last section, we use these results, together with an extension of \cite[VI.4, Proposition 1]{Ste}, to obtain the continuity of pseudodifferential operators with symbols in ${\rm S}^m_{\rho,\delta,N,N'}$ 
on mixed-norm Lebesgue spaces.

\smallskip

\noindent {\bf Notation.}
We use the following notation and definition for the Fourier transform of the function $u\in L^1(\R^d)$:
\[\hat u(\xi)={\cal F}u(\xi)=\int_{\R^d} e^{-ix\cdot \xi}u(x)d x\,.\]
When the area of integration is not specified, the whole space is assumed. We denote by ${\cal S}(\R^d)$ the Schwartz space of smooth rapidly decreasing functions on $\R^d$ and by ${\cal S}'(\R^d)$ its dual, the space of tempered distributions.
For $N,N'\in\N_0$ we use the following family of semi-norms on ${\cal S}(\mathbb R^d)$:
\[|\varphi|_{N,N'}= \sup_{|\alpha|\leq N,|\beta| \leq N'}\sup_{x \in \mathbb R^d}|x^{\alpha}\partial^{\beta}\varphi(x)|, \]
and by ${\cal S}_{N,N'}(\R^d)$ we denote the Banach space of all functions $\varphi\in C^{N'}(\R^d)$ for which $|\varphi|_{N,N'} < \infty$.
We denote Lebesgue spaces by $L^p(\R^d)$, locally Lebesgue spaces by $L^p_{loc}(\R^d)$ and Lebesgue spaces with compact support by $L^p_c(\R^d)$.
By $C$ we always denote a constant, even if it changes during calculation,
while $C_p$ is a constant depending on parameter $p$. By $\lfloor x \rfloor$ we denote the largest integer not greater than $x$, while $\lfloor x \rfloor_2$ is the largest even integer not greater than $x$. We also use the standard notation $m^+=\max\{m,0\}$.  
Finally, $|\cdot|_p$ is the standard $p$-norm on $\R^d$:
\[
|x|_p=\sqrt[p]{|x_1|^p+\cdots+|x_d|^p},\qquad |x|_\infty=\max\{|x_1|,\ldots,|x_d|\}\,,
\]
and we denote $|x|_2$ simply by $|x|$.

\section{The general framework}

Our main tool is the following theorem, proved in \cite[Theorem 1]{AIV}. To state the theorem the following notation is convenient:
Take $l\in\{0,...\, ,(d-1)\}$ and split $x=(\bar{x} ,x')=(x_1,\ldots,x_l;x_{l+1},\ldots,x_d)$.
Next define (with a slight abuse of notation, for simplicity) for $\bar{{\bf p}}=(p_1,\ldots,p_l)$
\[
L^{\bar{{\bf p}},\,p}(\R^d)=L^{(\bar{{\bf p}},p,\ldots,p)}(\R^d) \;\qquad \hbox{and} \qquad
	\|f\|_{\bar{{\bf p}},\,p}=\|f\|_{(\bar{{\bf p}},\,p,\ldots,\,p)}\;.
\]
Of course, for $l=0$ we take $\|f(\cdot,x')\|_{\bar{{\bf p}}}=|f(x')|$ and $\|f\|_{\bar{{\bf p}},\,p}=\|f\|_{L^p}.$
We also define (for each $l\in\{0,...\, ,(d-1)\},\,t>0$ and 
$y'\in\R^{d-l}$):
\[
{\cal F}_{l,t}^{y'}:=\Big\{f\in L^1_{loc}(\R^d):{\rm supp}\,f\subseteq \R^l\times\{x':|x'-y'|_\infty\leq t\}\;\;\&\;
\int_{\R^{d-l}} f(\bar{x},x') \, dx'	=	0\;({\rm ae}\;\bar{x}\in\R^l)\Big\}\,.
\]

\begin{theorem}
Assume that $A,A^\ast:L^\infty_c(\R^d)\to L^1_{loc}(\R^d)$ are formally adjoint linear operators, i.e.~such that
\[
(\forall\,\vp,\psi\in C^\infty_c(\R^d))	\quad
	\int_{\R^d}(A\vp)\overline\psi =\int_{\R^d}\vp\overline{A^\ast\psi}. 
\]
Furthermore, let us assume that (both for $T=A$ and $T=A^\ast$) there exist constants $N>1$ and $c_1>0$ satisfying
\begin{equation}\label{maincond}
(\forall\,l\in\{0,...\, ,(d-1)\})(\forall\,x'_0\in \R^{d-l})(\forall\,t>0)	\quad
	\int\limits_{|x'-x'_0|_\infty>Nt} \|Tf(\cdot,x')\|_{\bar{{\bf p}}} \,dx' \leq c_1\|f\|_{\bar{{\bf p}},1} \;,
\end{equation}
for any function $f\in L^\infty_c(\R^d)\cap{\cal F}_{l,t}^{x'_0}$ and any $\bar{{\bf p}}\in\langle1,\infty\rangle^l$.

If for some $q\in\langle1,\infty\rangle$ operator $A$ has a continuous extension to an operator from\/ $L^q(\R^d)$ to itself with norm $c_q$,
then $A$ can be extended by the continuity to an operator from\/ $L^{\bf p}(\R^d)$ to itself for any\/
${\bf p}\in \langle1,\infty\rangle^d$, with the norm
\begin{align}
\|A\|_{L^{\bf p} \to L^{\bf p}}
	& \leq \sum_{k=1}^d c^k \prod_{j=0}^{k-1} \max (p_{d-j},(p_{d-j}-1)^{-1/p_{d-j}}) (c_1+c_q) \nonumber \\
	& \leq c' \prod_{j=0}^{d-1} \max (p_{d-j},(p_{d-j}-1)^{-1/p_{d-j}}) (c_1+c_q), \nonumber
\end{align}
where $c$ and $c'$ are constants depending only on $N$ and $d$.
\end{theorem}

We also use \cite[Theorem 5.5]{IV} to conclude that for a pseudodifferential operator with symbol $\sigma\in {\rm S}^m_{\rho,\delta,N,N'}$, where $\delta\leq\rho$, $m\in[-d,0]$ and $N,N'\in2\N_0$ are such that
\begin{equation}\label{adreq}
N > \frac{(3-\delta)d+(5-\delta)(1-\delta)}{(1-\delta)^2},\quad N'>6d+12\,, 
\end{equation}
its formally adjoint operator exists and has symbol $\sigma^\ast\in {\rm S}^m_{\rho,\delta,M,M'}$, where $M,M'\in2\N_0$ satisfy
\begin{equation}\label{adreqq}
N-M>\frac{d+(\lfloor d\rfloor_2+2)\delta}{1-\delta},\quad N-M\geq\frac{-m+(1-\delta)d+(\lfloor d\rfloor_2+2)\delta}{1-\delta},\quad N'-M'\geq \lfloor d\rfloor_2+2\,.
\end{equation}
The claim is valid for all $m\leq0$ because it can easily be checked that the condition $m\geq -d$ used in \cite[Theorem 5.5]{IV} could be omitted.

\section{Boundedness of $\Psi$DO}

$S^m_{\rho,\delta,N,N'}$ is a nonsmooth variant of H\"ormander class $S^m_{\rho,\delta}$, $m\in\R$, $0\leq \rho \leq 1$, $0\leq \delta < 1$. It consists of all $\sigma:\R^d\times\R^d\to{\mathbb C}$ such that for all multi-indices $|\alpha|\leq N,|\beta|\leq N'$ it holds
\begin{equation}\label{srdn} 
(\forall x\in\R^d)(\forall \xi\in\R^d)\quad |\partial_x^{\alpha}\partial_{\xi}^{\beta}\sigma (x,\xi)| \leq C_{\alpha,\beta}\langle \xi \rangle^{m-\rho|\beta|+\delta|\alpha|}\,, 
\end{equation}
where $\langle \xi \rangle = (1+|\xi|^2)^{\frac12}$, $C_{\alpha,\beta}$ is a constant depending only on $\alpha$ and $\beta$ and where all partial derivatives are understood to be continuous.

For such symbols we denote the corresponding pseudodifferential operator $T_{\sigma}$ by
\begin{equation}\label{psido}
T_{\sigma}\varphi(x)=\int_{\R^d}e^{ix \cdot \xi}\sigma(x,\xi)\hat{\varphi}(\xi)\ \dbar \xi,\ \varphi \in {\cal S}(\R^d),
\end{equation}
where $\,\dbar\xi = (2\pi)^{-d}d\xi$. As $T_{\sigma}$ doesn't map ${\cal S}(\R^d)$ to ${\cal S}(\R^d)$ when symbols are not infinitely smooth, we cannot extend this definition to the space of tempered distributions, but we can extend it to mixed-norm Lebesgue spaces $L^{\bf p}(\R^d)$, ${\bf p}\in \langle1,\infty\rangle^d$ using a formula $\langle T_\sigma u , \varphi\rangle = \langle u , T_{\sigma^\ast}\varphi \rangle$, where $\langle \cdot , \cdot \rangle$ is the dual product, $u\in L^{\bf p}(\R^d)$, $\varphi\in{\cal S}(\R^d)$.
The only requirement for this to work is $T_{\sigma^\ast}\varphi\in L^{\bf p'}(\R^d)$. This requirement also guarantees that the continuous extension of $T_{\sigma}$ on $L^{\bf p}(\R^d)$, if it exists, should be given by the above duality formula. As $T_{\sigma^\ast}$ actually maps ${\cal S}(\R^d)$ to ${\cal S}_{M',M}(\R^d)$ \cite[Theorem 2.2]{IV}, a sufficient condition for the above requirement to be fulfilled is $M'\geq d$. This can be seen easily by writing
\[
\|T_{\sigma^\ast}\varphi\|_{L^{\bf p'}}=\Big(\int\langle x_d\rangle^{-1-\epsilon}\ldots\Big(\int\langle x_{1}\rangle^{-1-\epsilon} \Big(
\langle x_1\rangle^{\frac{1+\epsilon}{p_1}}\langle x_2\rangle^{\frac{1+\epsilon}{p_2}}\cdots\langle x_d\rangle^{\frac{1+\epsilon}{p_d}}|T_{\sigma^\ast}\varphi|\Big)^{p_1}\,dx_{1}\Big)^{\frac{p_2}{p_{1}}}\ldots\,dx_d\Big)^{\frac1{p_d}}
\]
and by noticing that we can take $\epsilon>0$ so small that $\frac{1+\epsilon}{p_i}\leq1$ for every $i\in\{1,2,\ldots,d\}$. This condition will be satisfied in our main result -- Theorem 3.2.

In this paper we prove that the continuous extension of $T_{\sigma}$ on $L^{\bf p}(\R^d)$ exists for $\rho>0$, $\delta\leq\rho$,
\begin{equation}\label{lppcond}
m\leq -(1-\rho)(d+1+\rho)\,,
\end{equation}
and for sufficiently smooth symbols. As the estimate (\ref{lppcond}) is more crude than (\ref{lpcond}) we expect that even better results are possible.

From the estimate (\ref{srdn}), it follows easily that for a fixed $x\in\R^d$ we have $\sigma(x,\cdot)\in{\cal S'}(\R^d)$ and so there is a $k(x,\cdot)\in{\cal S'}(\R^d)$ such that
$\widehat{k(x,\cdot)}=\sigma(x,\cdot)$. We call the tempered distribution $k(x,\cdot)$ a kernel of the operator $T_{\sigma}$ and using properties of the convolution and Fourier transform we can write (\ref{psido}) in the form
\begin{equation}\label{psidok}
T_{\sigma}\varphi(x)=k(x,\cdot)\ast\varphi\,.
\end{equation}
The kernel $k(x,\cdot)$ is a function away from the origin with the following estimates on its derivatives.

\begin{lemma}
Let $\sigma\in S^m_{\rho,\delta,N,N'}$, $\rho>0$. Then the kernel $k(x,z)$ satisfies
\begin{equation}\label{kest}
|\partial_x^{\alpha}\partial_{z}^{\beta}k (x,z)|\leq C_{\alpha,\beta,L}\cdot|z|^{-d-m-\delta|\alpha|-|\beta|-L}\,,\quad z\ne0\,,
\end{equation}
for all $|\alpha|\leq N$, $|\beta|\geq 0$ and
\begin{equation}\label{Lreq}
L\geq (1-\rho)\Big(\Big\lfloor\frac{d+m+\delta|\alpha|+|\beta|}{\rho}\Big\rfloor+1\Big)^+
\end{equation}
such that $N'\geq d+m+\delta|\alpha|+|\beta|+L>0$ and $N'>\frac{d+m+\delta|\alpha|+|\beta|}{\rho}$, and where $C_{\alpha,\beta,L}$
is a constant depending only on $\alpha,\beta$ and $L$.
\end{lemma}
\begin{proof}
We take a fixed nonnegative $\eta\in C^\infty_c(\R^d)$ such that $\eta(\xi)=1$ for $|\xi|\leq 1$ and $\eta(\xi)=0$ for $|\xi|\leq 2$. We also define $\zeta(\xi)=\eta(\xi)-\eta(2\xi)$.
Then we have the following decomposition (pointwise convergence in the sense of tempered distributions): 
\begin{equation*}
T_\sigma=\sum_{j=0}^\infty T_{\sigma_j}\,,
\end{equation*}
where $\sigma_0(x,\xi)=\sigma(x,\xi)\eta(\xi)$ and $\sigma_j(x,\xi)=\sigma(x,\xi)\zeta(2^{-j}\xi)$ for $j\geq 1$. Details of this so called dyadic decomposition are described in \cite[4.1--4.2]{Ste}. We denote kernels of the operators $T_{\sigma_j}\,(j\geq 0)$ by $k_j(x,\xi)$ and these are all smooth functions. Namely, as $\sigma_j(x,\xi)$ are functions with compact support they actually belong to ${\cal S}_{\infty,N'}(\R^d)$ and thus, by \cite[Lemma 2.1]{IV}, we have $k_j(x,\xi)\in{\cal S}_{N',\infty}(\R^d)\subseteq C^\infty(\R^d)$.

Obviously, $\sigma(x,\xi)=\sum_{j=0}^\infty\sigma_j(x,\xi)$ pointwise and so, using (\ref{srdn}) and Lebesgue dominated convergence theorem, we easily obtain that $k(x,\cdot)=\sum_{j=0}^\infty k_j(x,\cdot)$, with the sum converging for each fixed $x$ in the sense of tempered distributions. This, continuity of derivatives on the space of distributions and uniqueness of the limit in the space of distributions ensures that it is enough to prove that 
\begin{equation}\label{kjsum}
\sum_{j=0}^\infty |\partial_x^{\alpha}\partial_{z}^{\beta}k_j(x,z)|
\end{equation}
satisfies the estimate given by the right-hand side of (\ref{kest}).

So, we first find certain estimates on $|\partial_x^{\alpha}\partial_{z}^{\beta}k_j(x,z)|$. As $\sigma_j(x,\xi)$ have compact supports, we have
\[
k_j(x,z)=\int_{\R^d}\sigma_j(x,\xi)e^{i\xi \cdot z}\ \dbar \xi\,,
\]
and now for all multi-indices $|\alpha|\leq N$, $|\beta|\geq 0$ and $|\gamma|\leq N'$ it follows
\begin{equation}\label{kerjder}
(-iz)^\gamma\partial_x^{\alpha}\partial_{z}^{\beta}k_j(x,z)=\int_{\R^d}\partial_\xi^\gamma\Big((i\xi)^\beta\partial_x^{\alpha}\sigma_j(x,\xi)\Big)e^{i\xi \cdot z}\ \dbar \xi\,,
\end{equation}
and, using the fact that $\sigma_j$ is supported in $|\xi|\leq 2^{j+1}$ and for $j\ne 0$ also in $|\xi|\geq 2^{j-1}$, we can estimate the integrand above:
\begin{align}
\Big|\partial_\xi^\gamma\Big((i\xi)^\beta\partial_x^{\alpha}\sigma_j(x,\xi)\Big)\Big|&=\bigg|\sum_{\gamma'\leq\gamma}{\gamma \choose \gamma'}\partial_\xi^{\gamma'}(i\xi)^\beta\partial_x^{\alpha}\partial_\xi^{\gamma-\gamma'}\sigma_j(x,\xi)\bigg| \nonumber \\
&\leq C_\gamma\sum_{\gamma'\leq\gamma}|\xi|^{(|\beta|-|\gamma'|)_+}\bigg|\sum_{\gamma''\leq\gamma-\gamma'}{\gamma-\gamma' \choose \gamma''}\partial_x^{\alpha}\partial_\xi^{\gamma''}\sigma(x,\xi)\partial_\xi^{\gamma-\gamma'-\gamma''}\big(\zeta(2^{-j}\xi)\big)\bigg| \nonumber \\
&\leq C_{\alpha,\gamma}\sum_{\gamma'\leq\gamma}\sum_{\gamma''\leq\gamma-\gamma'}|\xi|^{(|\beta|-|\gamma'|)_+}\cdot\langle\xi\rangle^{m-\rho|\gamma''|+\delta|\alpha|}\cdot 2^{-j|\gamma-\gamma'-\gamma''|} \nonumber \\
&\leq C_{\alpha,\beta,\gamma}\sum_{\gamma'\leq\gamma}\sum_{\gamma''\leq\gamma-\gamma'}2^{j(|\beta|-|\gamma'|)}\cdot 2^{j(m-\rho|\gamma''|+\delta|\alpha|)}\cdot 2^{-j|\gamma-\gamma'-\gamma''|} \nonumber \\
&=C_{\alpha,\beta,\gamma}\sum_{\gamma'\leq\gamma}\sum_{\gamma''\leq\gamma-\gamma'}2^{j(m+\delta|\alpha|+|\beta|-\rho|\gamma|+(\rho-1)(|\gamma|-|\gamma''|)} \nonumber \\
&\leq C_{\alpha,\beta,\gamma}\cdot 2^{j(m+\delta|\alpha|+|\beta|-\rho|\gamma|}\,, \nonumber
\end{align}
where we have used $\rho\leq 1$ in the last estimate. From (\ref{kerjder}) we now obtain
\begin{align}
|z^\gamma\partial_x^{\alpha}\partial_{z}^{\beta}k_j(x,z)|&\leq C_{\alpha,\beta,\gamma}\cdot 2^{j(m+\delta|\alpha|+|\beta|-\rho|\gamma|} \cdot 2^{jd} \nonumber \\
&=C_{\alpha,\beta,\gamma}\cdot 2^{j(d+m+\delta|\alpha|+|\beta|-\rho|\gamma|}\,, \nonumber
\end{align}
and by taking supremum over $|\gamma|=M$ (for instance, we can take $\gamma$ such that $|z^\gamma|=(\max_{1\leq i\leq d}|z_i|)^M$) 
we finally get (for $|\alpha|\leq N$, $|\beta|\geq 0$ and $\N_0\ni M \leq N'$)
\begin{equation}\label{kjest}
|\partial_x^{\alpha}\partial_{z}^{\beta}k_j(x,z)|\leq C_{\alpha,\beta,M}\cdot |z|^{-M}\cdot 2^{j(d+m+\delta|\alpha|+|\beta|-\rho M)}\,.
\end{equation}
In order to estimate (\ref{kjsum}) we first consider the case $0<|z|\leq 1$. We split (\ref{kjsum}) into two parts: 
\[
S_1=\sum_{2^j\leq|z|^{-1}}|\partial_x^{\alpha}\partial_{z}^{\beta}k_j(x,z)|\qquad\hbox{and}\qquad S_2=\sum_{2^j>|z|^{-1}}|\partial_x^{\alpha}\partial_{z}^{\beta}k_j(x,z)|\,.
\]
Because of (\ref{kjest}) (for $M=0$) and $j\leq\log_2|z|^{-1}$ we have
\[
S_1\leq C_{\alpha,\beta}\sum_{2^j\leq|z|^{-1}} 2^{j(d+m+\delta|\alpha|+|\beta|}\leq C_{\alpha,\beta}\cdot
\begin{cases}
|z|^{-d-m-\delta|\alpha|-|\beta|}, & \text{if }\;d+m+\delta|\alpha|+|\beta|>0 \\ 
1+\ln(|z|^{-1}), & \text{if }\;d+m+\delta|\alpha|+|\beta|\leq 0\,.
\end{cases}
\]
Using an elementary inequality $\ln(1+x)\leq\frac1\alpha x^\alpha$, valid for $\alpha\in\langle0,1]$ and $x\geq0$, we finally obtain
\[
S_1\leq C_{\alpha,\beta,L}\cdot |z|^{-d-m-\delta|\alpha|-|\beta|-L}\,,
\]
for all $L\geq0$ such that $d+m+\delta|\alpha|+|\beta|+L>0$. To estimate $S_2$ we use (\ref{kjest}) with $M>\frac{d+m+\delta|\alpha|+|\beta|}{\rho}$:
\begin{align}
S_2&\leq C_{\alpha,\beta,M}\cdot |z|^{-M}\cdot \sum_{2^j>|z|^{-1}}2^{j(d+m+\delta|\alpha|+|\beta|-\rho M)} \nonumber \\
&\leq C_{\alpha,\beta,M}\cdot |z|^{-M-d-m-\delta|\alpha|-|\beta|+\rho M} \nonumber \\
&\leq C_{\alpha,\beta,L}\cdot |z|^{-d-m-\delta|\alpha|-|\beta|-L}\,, \nonumber
\end{align}
for $L\geq(1-\rho)M$, which is exactly the requirement (\ref{Lreq}) and the proof is complete in the case $|z|\leq1$.

To estimate (\ref{kjsum}) in the case $|z|>1$ we again use (\ref{kjest}) with $M>\frac{d+m+\delta|\alpha|+|\beta|}{\rho}$ to get
(if we also take $M\geq d+m+\delta|\alpha|+|\beta|+L$, which is possible under requirements of the lemma)
\[
\sum_{j=0}^\infty |\partial_x^{\alpha}\partial_{z}^{\beta}k_j(x,z)|\leq C_{\alpha,\beta,M}\cdot |z|^{-M}\leq C_{\alpha,\beta,L}\cdot |z|^{-d-m-\delta|\alpha|-|\beta|-L}\,.
\]
\end{proof}

If $N'>\Big(\frac{d+m}\rho\Big)^+$, the bound from the previous lemma provides us with the following integral representation
\begin{equation}\label{Tintrep}
T_{\sigma}f(x)=\int_{\R^d} k(x,x-y)f(y) \, dy\,,
\end{equation}
for any $f\in C_c^\infty(\R^d)$ and $x\not\in {\rm supp} f$. Obviously, we can also take $f\in L_c^\infty(\R^d)$ and the representation is than valid for a.e.~$x\not\in {\rm supp} f$.

Moreover, because of the $L^2$ continuity result by Coifman and Meyer mentioned in the Introduction and density of $C_c^\infty(\R^d)$ in $L^2(\R^d)$, we can easily show that, for operators of order zero and sufficiently large $N'$,
this representation remains valid for any
$f\in L^2(\R^d)$ and a.e.~$x\not\in {\rm supp} f$.

We are now in a position to prove the main theorem of this paper.
\begin{theorem}
Let $\sigma\in S^m_{\rho,\delta,N,N'}$, $[0,1\rangle\ni\delta\leq\rho\in\langle0,1]$ and $m\leq -(1-\rho)(d+1+\rho)$. If 
\[
N > \frac{(3-\delta)d+(5-\delta)(1-\delta)}{(1-\delta)^2},\quad N'>6d+12\,,
\]
then $T_\sigma$ is bounded on $L^{\bf p}(\R^d)$, ${\bf p}\in \langle1,\infty\rangle^d$.
\end{theorem}
\begin{proof}
Because of Theorem 2.1 and known $L^2$ continuity result, we only need to show that $T_\sigma$ and $T_{\sigma^\ast}$ satisfy
the estimate (\ref{maincond}). Moreover, because of (\ref{adreq})-(\ref{adreqq}) it is enough to show that $T_\sigma$ satisfies the estimate (\ref{maincond}) -- the result for $T_{\sigma^\ast}$ is then an easy consequence.

We now follow the proof of \cite[Lemma 3]{AIV}. For an arbitrary $N>1$, using (\ref{Tintrep}), two linear changes of variables, and the Minkowski inequality, we get
\begin{align}
&\!\!\!\!\!\!\!\!
\int\limits_{|x'-x'_0|_{\infty}>Nt}\!\!\|T_\sigma f(\cdot,x')\|_{\bar{\bf p}} \, dx' \nonumber \\
&=\!\!\!\!
\int\limits_{|x'-x'_0|_{\infty}>Nt}\!\!\Big\|\int\int k(\cdot,x',\cdot-{\bar y},x'-y')f({\bar y},y') \, d{\bar y} \, dy'\Big\|_{\bar{\bf p}} \, dx' \nonumber \\
&=\!\!\!\!
\int\limits_{|x'|_{\infty}>Nt}\!\!\!\!\Big\|\int\int\Big(k(\cdot,x'\!+\!x'_0,\cdot-\!{\bar y},x'\!-\!y')-k(\cdot,x'\!+\!x'_0,\cdot-\!{\bar y},x')\Big)
    f({\bar y},y'\!+\!x'_0) \, d{\bar y} \, dy'\Big\|_{\bar{\bf p}} \, dx' \nonumber \\
&=\!\!\!\!
\int\limits_{|x'|_{\infty}>Nt}\!\!\!\!\Big\|\int\int\Big(k(\cdot,x'\!+\!x'_0,{\bar y},x'\!-\!y')-k(\cdot,x'\!+\!x'_0,{\bar y},x')\Big)
    f(\cdot-\!{\bar y},y'\!+\!x'_0) \, d{\bar y} \, dy'\Big\|_{\bar{\bf p}} \, dx' \nonumber \\
&\leq\!\!\!\!
\int\limits_{|x'|_{\infty}>Nt}\!\!\!\!\int\int\Big\|\Big(k(\cdot,x'\!+\!x'_0,{\bar y},x'\!-\!y')-k(\cdot,x'\!+\!x'_0,{\bar y},x')\Big)
    f(\cdot-\!{\bar y},y'\!+\!x'_0)\Big\|_{\bar{\bf p}} \, d{\bar y} \, dy' dx'=I. \nonumber
\end{align}
Furthermore, using (\ref{kest}) and the Mean value theorem, for $|x'|_{\infty}>Nt\,$ and $\,|y'|_{\infty}\leq t$ we have an estimate
\begin{align}
|k(\bar x,x'+x'_0,\bar y,x'-y')-k(\bar x,x'+x'_0,\bar y,x')|&=
    |\nabla_{y'}k(\bar x,x'+x'_0,\bar y,x'-\vartheta y')\cdot y'| \nonumber \\
&\leq C|(\bar y,x'-\vartheta y')|^{-d-m-1-L}\,|y'|_\infty \nonumber \\
&\leq C|(\bar y,x'-\vartheta y')|^{-d-m-1-L}\,t, \nonumber
\end{align}
for some constants $\vartheta\in\langle0,1\rangle$, $C>0$ and $L$ satisfying assumptions of Lemma 3.1 (the optimal value of $L$ will be determined later).
Now, using the assumption on the support of function $f$ we continue the estimate
\begin{align}
I&\leq Ct\int\limits_{|x'|_{\infty}>Nt}\;\int\limits_{|y'|_{\infty}\leq t}\int|(\bar y,x'-\vartheta y')|^{-d-m-1-L}\,\|f(\cdot-\bar y,y'+x'_0)\|_{\bar {\bf p}} \, d\bar y \, dy' dx' \nonumber \\
&= Ct\int\limits_{|x'|_{\infty}>Nt}\;\int\limits_{|y'|_{\infty}\leq t}\int|(\bar y,x'-\vartheta y')|^{-d-m-1-L}\,\|f(\cdot,y'+x'_0)\|_{\bar {\bf p}} \, d\bar y \, dy' dx' \nonumber \\
&= Ct\int\limits_{|y'|_{\infty}\leq t}\|f(\cdot,y'+x'_0)\|_{\bar {\bf p}}\int\limits_{|x'|_{\infty}>Nt}\int|(\bar y,x'-\vartheta y')|^{-d-m-1-L} \, d\bar y \, dx' dy'. \nonumber
\end{align}
To conclude the proof we need to check that
\[
I\!I=t\int\limits_{|x'|_{\infty}>Nt}\int|(\bar y,x'-\vartheta y')|^{-d-m-1-L} \, d\bar y \, dx'
\]
is bounded on $|y'|_{\infty}\leq t,$ for arbitrary $N>d$. Indeed, we have
\begin{align}
d^{-d-1}I\!I&\leq t\int\limits_{|x'|>Nt}\int(|\bar y|_1+|x'-\vartheta y'|_1)^{-d-m-1-L} \, d\bar y \, dx' \nonumber \\
&\leq t\int\limits_{|x'|>Nt}\int(|\bar y|_1+|x'|_1-|y'|_1)^{-d-m-1-L} \, d\bar y \, dx' \nonumber \\
&\leq t\int\limits_{|x'|>Nt}\int(|\bar y|_1+|x'|_1-dt)^{-d-m-1-L} \, d\bar y \, dx' \nonumber \\
&= t^{-m-L}\int\limits_{|x'|>N}\int(|\bar y|_1+|x'|_1-d)^{-d-m-1-L} \, d\bar y \, dx' \,, \nonumber
\end{align}
where in the last step we used a linear change of variables. It is now obvious that this can be bounded (uniformly for $t>0$) only in the case $L=-m\,$ 
and in that case the above integral is indeed finite.

Thus, according to Lemma 3.1 we need $N'\geq d+1$, $N'>\frac{d+m+1}{\rho}$ and
\[
-m\geq (1-\rho)\Big(\Big\lfloor\frac{d+m+1}{\rho}\Big\rfloor+1\Big)^+\,.
\]
The last requirement is obviously valid in the case $m<-d-1$. Otherwise, it is sufficient to have
\[
-m\geq (1-\rho)\Big(\frac{d+m+1}{\rho}+1\Big)\,,
\]
which is equivalent to $m\leq-(1-\rho)(d+1+\rho)$ -- the assumption stated in the theorem.

Of course, the same estimate should be valid also for $T_{\sigma^\ast}$. Because of (\ref{adreq})-(\ref{adreqq}) this means that we also need $M'\geq d+1$ and $M'>\frac{d+m+1}{\rho}$, where $M,M'$ are defined there. In the end, we notice that it is enough to take $M=0$ and that the assumptions of this theorem are sufficient in that case.

\end{proof}

\end{document}